\documentclass[11pt]{amsart} 
\usepackage{amssymb,amsmath,latexsym,enumerate,graphicx,bbm,mathptmx,microtype,cite,xcolor,ifthen}
\usepackage{tabularx}

\allowdisplaybreaks

\newtheorem{theorem}{Theorem} 

\newtheorem{exam}{Example}

\newtheorem*{rem}{Remarks}

\theoremstyle{definition}
\newtheorem{definition}[theorem]{Definition}

\newcommand\commentout[1]{}

\newcommand\lcm{\operatorname{lcm}}

\newcommand\ZZ{\mathbb{Z}}
\newcommand\QQ{\mathbb{Q}}

\newcommand\rat{\!\! : \!\!}
\newcommand\sh{$\sharp$}
\newcommand\fl{$\flat$}

\makeatletter 
\newtheorem*{rep@theorem}{\rep@title}\newcommand{\newreptheorem}[2]{%
\newenvironment{rep#1}[1]{%
\def\rep@title{\bf #2 \ref{##1}}%
\begin{rep@theorem}}%
{\end{rep@theorem}}}
\makeatother
\newreptheorem{theorem}{Theorem}

\newcounter{teach}
\begin{document}

\title{Musical Chords by the Numbers}  

\author{Matthias Beck}
\author{Emily Clader}
\address{Department of Mathematics\\
         San Francisco State University\\
         San Francisco, CA 94132\\
         U.S.A.}
\email{[mattbeck,eclader]@sfsu.com}


\begin{abstract}
The mathematics of musical intervals and scales has been extensively studied. Vastly simplified, our ears
seem to prefer intervals whose frequency ratios have small numerator and denominator, such as 2:1
(octave), 3:2 (perfect fifth), 4:3 (perfect fourth), and so on. While there also have been
numerous studies on the mathematics of musical chords, very few of them consider a model that measures
consonance/dissonance of a given chord in analogy with this simple-fractions perspective. Our aim is to
develop a measure for the consonance of a chord with crucial symmetry features, including invariance under
chord translation, inversion, and interval sets. We apply our model to chords in various musical scales and compare it
to existing models and empirical studies.
\end{abstract}



\date{19 September 2025}
\thanks{We thank Federico Ardila, James T.\ Smith, Frieder Stolzenburg, and two anonymous referees for helpful comments and conversations about this project.}

\maketitle


\section{Overture}

The mathematics of musical scales---and, by extension, intervals between two notes---has been
extensively studied; see, e.g., the classic~\cite{bohlen,huron} or the more recent~\cite{cladertwelvetones,schettlerxenharmonic}.  The notes
of a scale can naturally be recorded numerically via their frequencies, and an interval is then encoded via the frequency ratio between its two notes. Vastly simplified (and one has to be careful, for various reasons, some of which we summarize below), our ears
seem to prefer ratios with small numerator and denominator, such as $2 \rat 1$ (octave), $3 \rat 2$ (perfect
fifth), $4 \rat 3$ (perfect fourth), and so on. Some relevant background on the mathematics of intervals can be
found in, e.g.,~\cite{bensonmusic,walkerdon}.

\begin{figure}[ht] 
  \includegraphics[height=.57in]{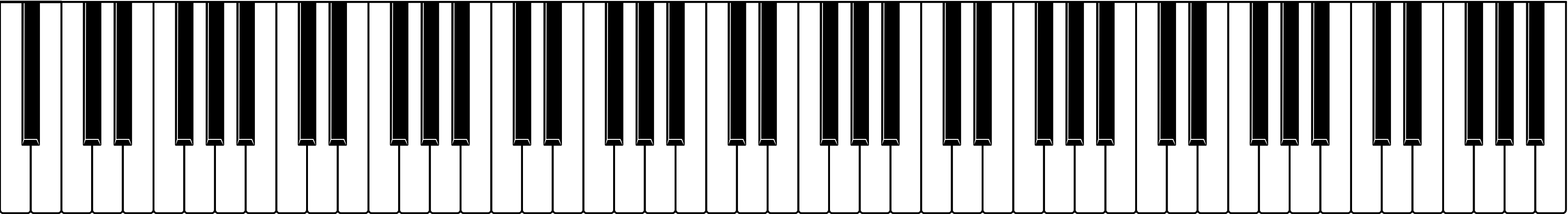}
\caption{Our favorite depiction of $\ZZ_{ 12 }$. \tiny [Figure credit: Wikimedia Commons,
Piano-full-en.svg]}\label{fig:pianokeys}
\end{figure}

Generalizing beyond intervals, there is also much literature on the mathematics of musical chords consisting of three or more notes.  This literature gives interesting group-theoretic connections, both in music theory (see, e.g.,~\cite{morristwelvetone}) and algebra
(see, e.g.,~\cite{cransfioresatyendra}). These tend to concern transformations that a chord can undergo, such as transpositions and inversions. It is a short step to model, say, a twelve-note scale via
$\ZZ_{ 12 }$ (see Figure~\ref{fig:pianokeys}), and
so these mathematical aspects of chords are independent of the actual frequency ratios of the scale.

But why do we perceive, say, a major triad as more pleasant than a chord consisting of three notes
separated by two half tones? One could make an argument that the appreciation of chords is to a large
extent culturally based. On the other hand, acoustics (see, e.g.,~\cite{hallacoustics}) offers
some explanations, for instance by looking at the wave forms of chords containing a perfect fourth. (Many piano
tuners start with fourths for the same reason, though they listen for beats in order to
achieve an
equal-tempered scale.) This gives rise to much interesting mathematics (see, e.g.,~\cite{cubarsi}) and
statistics (see, e.g.,~\cite{eerolalahdelma,masinalopresti}), not to mention empirical (see,
e.g.,~\cite{johnsonlairdkangleong,robertsconsonant}) and psychophysical (see,
e.g.,~\cite{cookfujisawa,parncutt,terhardtstollseewann}) results on how consonant\footnote{
Consonance/dissonance of musical chords are not mathematically well-defined concepts---the composer and music theorist Paul Hindemith famously wrote ``the two concepts have never been
completely explained, and for a thousand years the definitions have varied.''
We will resist the urge to put the words \emph{consonant} and \emph{dissonant} between quotation
marks; nevertheless, we use them only in terms of (vaguely) comparing two given chords.
} 
different chords are perceived to be.
However, the mathemusical literature seems to contain few studies on musical chords in analogy with the
aforementioned simple-fractions perspective on intervals, and our aim is to fill in some of these gaps.
The two places in the literature in which we could find existing work along these lines\footnote{
We should mention that Birkhoff~\cite[Chapter V]{birkhoffaestheticmeasure} conceptualized a notion of \emph{aesthetic measure} of a chord, though it is arguably much less refined.
} 
are Brefeld's website on mathematics in our
everyday world \cite{brefeld} and Stolzenburg's papers on harmony perception
\cite{stolzenburgproc,stolzenburgjmm}; we will review their work in Section~\ref{sec:review}.

The aim of our work is to develop a new model by which a numerical value can be assigned to any chord
that, in some sense, measures how consonant the chord is.  Our model will share certain
features with Stolzenburg's preexisting notion of \emph{harmonicity}, but in contrast to
Stolzenburg's work, the numerical value that we assign to a chord will be unchanged if the chord
is translated or inverted, and for this reason, we refer to our concept as {\it symmetric
harmonicity}.  This notion will rely on first fixing a scale---that is, on specifying the
numerical frequency that will correspond to each note---and different scales will yield different
numerics.  After developing our model, we will test it out on different scales and compare our findings to those in the literature, arguing that the notion of symmetric harmonicity closely mirrors an intuitive notion of consonance when the scale in question has roughly equal frequency ratios between each pair of consecutive notes.

A few words on our assumptions are in order.  First, because the notes produced by actual musical instruments
are not perfect sine waves but rather sums of sine waves, they are more accurately encoded not by a single
frequency but by a fundamental frequency along with a series of overtones; the frequencies and relative
amplitudes (loudness) of those overtones are, at least roughly, what produces the unique sound or \emph{timbre} of each instrument.  In the seminal book \cite{sethares}, Sethares has argued that consonance depends not only on the pitch (frequency) of notes but on their timbre, so preferences about scales and chords should be chosen based on the timbre of the particular instrument at hand.  Depending on the timbre, even octaves need not have the special ``most consonant'' status that they traditionally hold; see, for example, the work of Bohlen, Mathews, Roberts, and Pierce \cite{mathewsrobertspierce, mathewspierce}, which considers a scale where the role of the octave is instead played by a 3:1 frequency ratio.

We will ignore these subtleties in our model, assuming throughout that the timbres at hand are ``harmonic
sounds'' \cite[p.~3]{sethares}, meaning their overtones series consists of integer multiples of the fundamental tone---which is true of the guitar, piano, and many other Western instruments---and therefore, if a note is identified with its fundamental frequency, the simple-fractions perspective on intervals outlined in the first paragraph indeed gives a rough measure of perceived consonance.  Furthermore, we will assume in our model that the scales at hand are \emph{just scales}, meaning scales in which all frequency ratios are rational numbers.  One side effect of working with just scales is that, though we generally prefer scales where the frequency ratios between successive notes are roughly equal, they can never be precisely equal; a scale
with the latter property is called \emph{equal-tempered}, and if it has twelve notes, the unique frequency ratio is the irrational number $\sqrt[12]{2}$.  (See \cite{bolandhughston} for a mathematical interplay
between just and equal-temperament scales.)  Lastly, we assume that our scales have twelve notes, though this is merely for convenience, and our model can easily be adapted to just scales with any number of notes.

Like most (all?) mathematical texts on musical subjects, we stress that we have no intention to
capture, much less understand, the way different chords feel; this is (way) beyond any arithmetic
nature, not to mention the rich music theory and accompanying music cognition that gets much closer to this issue (see, e.g., \cite{huronbook}). Along similar
lines, we also clarify from the start that our model does not constitute any absolute measure: some 
features will seem arbitrary, and we feel that this is fine.  We are mathematicians studying a mathematical structure that we find intriguing, and we invite the reader to come up with their own
interpretation of what this mathematics might mean music-theoretically.


\section{Prelude: Previous Studies}\label{sec:review}

We now briefly review two previous \emph{ans\"atze} for assigning a numerical value to a given chord.
Like ours, they depend on a given (just) scale, and for illustration's purpose, it will be
helpful to have a sample scale to follow along with.  The sample scale we will use is due to Johannes Kepler and is known as \emph{Kepler's Monochord No.~2} (transposed down a fifth) \cite{barbour}.  In Table \ref{tab:classicjustscale}, we describe this scale by indicating, for each note, its frequency ratio to the bottom note of the scale (though it is worth noting that the frequency ratios are independent of the chosen bottom note, so the labels in the first row of the table are for illustration only).\footnote{\label{justfootnote}
Kepler's just scale is constructed to satisfy two conditions.  First, the three main major chords
(C-E-G, F-A-C, and G-B-D) all appear with frequency ratio $4 \rat 5 \rat 6$; this is precisely
the ratio in which the major chord appears within the harmonic spectrum of a given tone.  Second,
the scale is symmetric, in the sense that the frequency ratio between the bottom C and the note
$s$ semitones above it is equal to the frequency ratio between the top C and the note $s$
semitones below it.  These two conditions force the frequency (relative to the bottom C) of all
the notes in the scale except the tritone (F\sh), whose frequency ratio is computed as the product of the ratios for a major third and a whole tone.}
%

\begin{table}[ht]
\begin{center}
\begin{tabular}{c|c|c|c|c|c|c|c|c|c|c|c}
C & C\sh & D & E\fl & E & F & F\sh & G & G\sh & A & B\fl & B \\
\hline 
$1 \rat 1$ & $16 \rat 15$ & $9 \rat 8$ & $6 \rat 5$ & $5 \rat 4$ & $4 \rat 3$ & $45 \rat 32$ & $3
\rat 2$ & $8 \rat 5$ & $5 \rat 3$ & $16 \rat 9$ & $15 \rat 8$
\end{tabular}
\end{center}
\vspace{10pt}
\caption{The frequency ratios of Kepler's just scale.}\label{tab:classicjustscale}
\end{table}





Given such a scale, Brefeld~\cite{brefeld} computes, for each interval, what he calls its \emph{consonance value}, defined
as the geometric mean of the numerator and denominator of the respective frequency ratio. Thus,
e.g., a perfect fifth in Kepler's just scale has consonance value $\sqrt 6$.  The (over-simplified) idea behind this value is that a smaller consonance value indicates a ``simpler'' frequency ratio, and tones with simple frequency ratios share more frequencies in their harmonic spectra, leading them to generally sound more consonant together.\footnote{
Arguing about consonance/dissonance via the agreement/disagreement of frequencies in the harmonic
spectra goes back to at least Helmholtz~\cite{helmholtz}, yielding the notion of \emph{roughness}.
Naturally, it does not explain perceived consonance of, say, two sine waves, but more recent work
studied neuronal periodicity detection~\cite{ebeling}, which serves as a point of departure for
Stolzenburg's work.
}

Brefeld then extends this concept to the \emph{consonance value} of a chord, defined as the geometric
mean of the consonance values of all intervals appearing in the chord.  For example, the major triad C-E-G contains three intervals: a major third (C-E), a minor third (E-G), and a perfect fifth (C-G), so its consonance value with respect to Kepler's just scale equals
\[
  \sqrt[3]{\sqrt{ 5 \cdot 4} \cdot  \sqrt{ 6 \cdot 5} \cdot  \sqrt{ 3 \cdot 2}} \ \approx \ 3.91 \, .
\]
Again, a smaller consonance value suggests that the chord is more consonant.
We note that if we simply want to compare chords, then the roots are irrelevant, and we may instead consider a modified consonance value by simply
multiplying all numerators and denominators of all intervals appearing in a given chord.

Stolzenburg~\cite{stolzenburgproc,stolzenburgjmm} developed the following alternative. He computes
what he calls the \emph{harmonicity} or the \emph{relative periodicity} of a chord based on the frequency ratios of its components
relative to the lowest tone in the chord. If the chord contains $k$ tones whose frequency ratios to the lowest tone are 
$\frac{ a_1 }{ b_1 } = 1$, $\frac{ a_2 }{ b_2 }$, \dots, $\frac{ a_k }{ b_k }$, then the relative periodicity
of the chord equals $\lcm(b_1, b_2, \dots, b_k)$.
The (again simplified) idea behind this value is that it equals the period length of the
superposition of the sinusoids corresponding to the components of the chord, measured relative to the
period length of the lowest tone.  As an example, the major triad C-E-G in Kepler's just scale has relative frequencies $\frac{1}{1}, \; \frac{5}{4}, \; \frac{3}{2}$, so its harmonicity equals
\[\lcm(1,4,2) = 4.\]
Similarly to Brefeld's measure, a smaller relative periodicity suggests that the chord is more consonant.

\begin{figure}[ht] 
  \includegraphics[height=.52in]{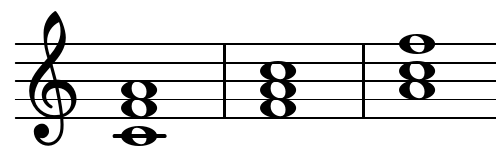}
\caption{F-major inversions.}\label{fig:inversions}
\end{figure}

Both Brefeld's and Stolzenburg's results correlate well with empirical data on harmony perception.
Stolzenburg has extended his work~\cite{stolzenburgjmm} into cognitive science, applying results from psychophysics and neuroacoustics.
We will present another alternative model on the mathematical side, one that emphasizes that certain chords are generally perceived to be ``equivalent'' to one another. 
To give an indication, we note that both Brefeld and Stolzenburg compute different values for inversions of the same chord (e.g., in Figure~\ref{fig:inversions}).
(Stolzenburg accounts for this by taking the average of the relative periodicities over all
inversions.)  Furthermore, in both Brefeld's and Stolzenburg's models, translating a chord generally changes its value (e.g., a C-major and G-major chord would have different values); this is simply because both models work directly with the frequency ratios between notes in the chord, and since the scale is not equal tempered, these frequency ratios are not translation invariant.  Mathematically (and also musically), we are seeking a model where the consonance measure is invariant under both inversions and translations. We will introduce this and other desiderata for our model next. 


\section{Bridge: Equivalence Assumptions}


Similarly to Brefeld's consonance value and Stolzenburg's harmonicity, our model will take as
input a fixed (just) scale, and will assign to each triad in that scale a numerical value that we
call its \emph{symmetric harmonicity}.  The key feature that distinguishes it from previous
models is that we will ask that symmetric harmonicity respects the following three notions of \emph{equivalence} of triads:
\begin{enumerate}
\item[(a)] chords that are translates of each other (e.g., C-minor and D-minor, obtained by translating by two scale steps) are equivalent;
\item[(b)] chords that are inversions of each other (e.g., in Figure~\ref{fig:inversions})
are equivalent;
\item[(c)] chords with the same set of pairwise intervals (e.g., C-minor and C-major, each of which
consists of a major third, a minor third, and a perfect fifth) are equivalent.
\end{enumerate}

Musically speaking, all of these are somewhat natural, though perhaps it is worth acknowledging their drawbacks at the outset.  First, since our scales are not equal-tempered, the actual frequency ratios between the notes in translated chords (such as between the notes in a C-minor versus a D-minor chord) will be different from one another, but we are insisting that their symmetric harmonicity values be the same.  In this sense, one might view the most natural scales for our model as being just scales that are ``close'' to equal-tempered; we will mention a particularly nice example of such a scale below.

Second, it should be acknowledged that consonance as experienced by actual humans is not translation invariant: the same chord generally sounds more dissonant when played at a lower frequency.  This can be explained by the fact that the critical bands of frequencies---which control the range of frequencies our ears discern as the ``same pitch''---get wider as the frequencies get higher; see \cite[page~44]{sethares} or \cite{plomp, zwickerflottorpstevens}.  To make our mathematical assumptions musically reasonable, then, we should assume that all pitches under consideration are within a fairly limited range of frequencies.

Finally, we note that equivalence (c) leads to the somewhat infamous question in harmony theory of why major and minor triads, despite being comprised of the same sets of intervals, sound different in character to our ears.  In spite of this difference, our model will view such triads as equally consonant.\footnote{Brefeld's consonance value also does not distinguish between major and minor triads, though Stolzenburg's harmonicity---since it considers not all pairwise intervals within a triad but only the intervals to the base tone---does distinguish between them.}  We view this equivalence as appealing mathematically, and perhaps one could take the results of this paper as an indication that musically, such an equivalence is admittedly imperfect but not entirely unreasonable.

Taking assumptions (a)--(c) as given, let us now represent them mathematically.  As mentioned in the introduction, one might view the space of scale-steps (prior to fixing a particular scale) as the set $\ZZ_{12} = \{0,1,2,\ldots, 11\}$---say by identifying C with 0, C\sh~with 1, D with 2, and so on\footnote{
Musicians think of C, C$\sharp$, etc., as a labeling of the \emph{chroma} of a given tone, and
organizing those by octaves yields a \emph{pitch class}---a wonderful example of the mathematical
concept of an equivalence class.
}---and from here, the space $T$ of triads is
\[T = \{(x_0, x_1, x_2) \; | \; x_0 < x_1 < x_2\} \subseteq \ZZ_{12}^3.\]
(Alternatively, one can view $T$ as the result of removing the diagonals $x_0 = x_1$, $x_0 = x_2$, and $x_1 =
x_2$ from $\ZZ_{12}^3$ and then taking a quotient by the symmetric group $S_3$; from this perspective, one sees
that $T$ embeds into $\ZZ_{12}^3/S_3$, which is a discrete version of an orbifold\cite{tymoczko}.)  Since
translation equivalence allows one to assume that $x_0 = 0$, the space of triads up to equivalence (a) can be identified with the set
\[T_0 = \{(x_1, x_2) \; | \; 0 < x_1 < x_2\} \subseteq \ZZ_{12}^2.\]

\begin{figure}[ht]
  \includegraphics[width=6.41in]{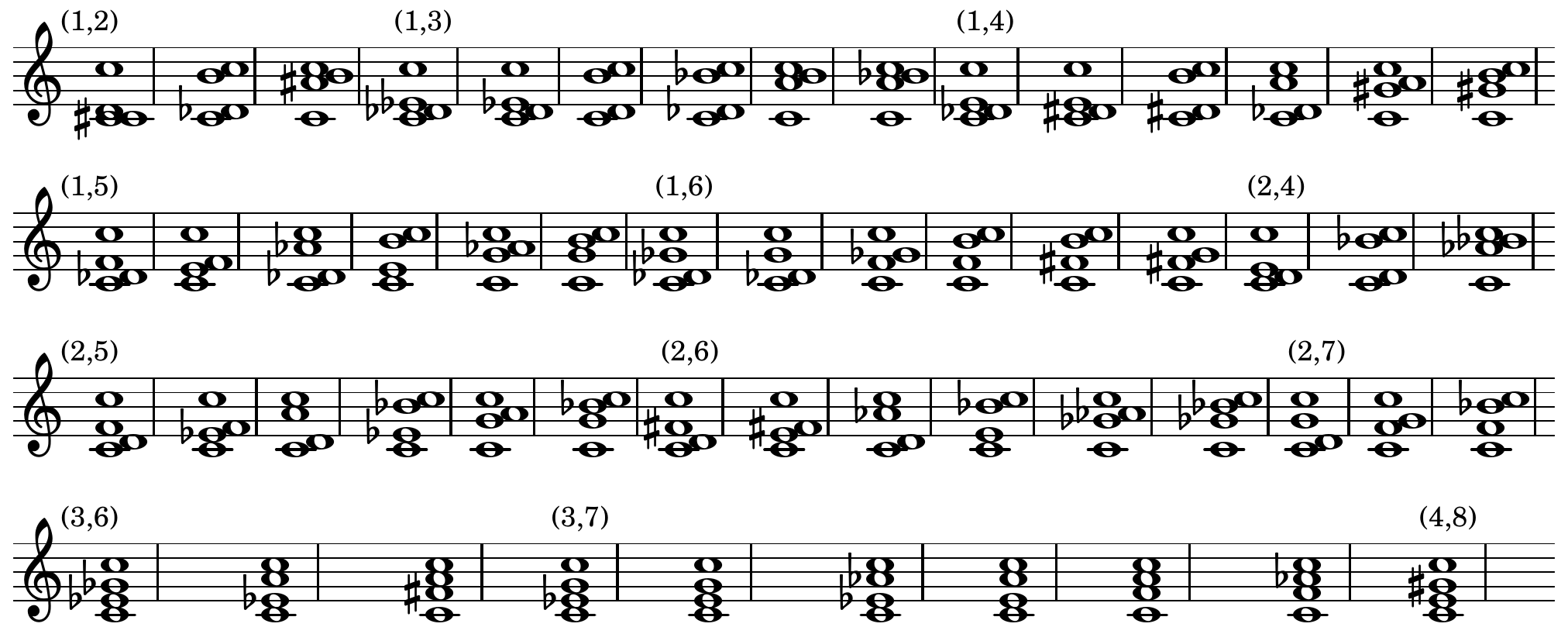}
\caption{All equivalence classes of triads.  For each equivalence class, we include all of its representatives $(x_1, x_2) \in T_0$, and we label it by one such representative.}\label{fig:allequclasses}
\end{figure}

Considering all three types of equivalence, it becomes an amusing combinatorial problem to count the number of equivalence classes of triads.  Viewing a triad, up to the translation-equivalence (a), as a triple $(0,x_1, x_2) \in T$, we see that the equivalences (b) and (c) leave unchanged the set
\[\{x_1-0, \; x_2-x_1, \; 0 -x_2\} \subseteq \ZZ_{12}\]
of intervals between consecutive notes (where we double the bottom note C at the octave), so we may organize the equivalence classes according to the number of elements in this set.  There are seven equivalence classes 
consisting of triads with three distinct intervals (each with six representatives in $T_0$), four equivalence classes 
consisting of triads with two
distinct intervals (each with three representatives in $T_0$), and one class 
in which all the notes are equally-spaced (with a single representative in $T_0$); 
Figure~\ref{fig:allequclasses} shows the complete set of equivalence classes (where we show
the bottom note C doubled at the octave).

%


\section{Exposition: Symmetric Harmonicity}

The way we will assign symmetric harmonicities to triads is to first assign them to intervals, and then to define the symmetric harmonicity of a chord as the product of the symmetric harmonicities of all intervals within it.

For intervals, in light of the translation-invariance (a), we need only assign symmetric harmonicities to intervals with bottom note C.  We will denote the symmetric harmonicity between C and the note $n$ scale-steps above it by $h_n$, so that $h_1$, for example, denotes the symmetric harmonicity of a semitione (the interval between C and C\sh), $h_2$ denotes the symmetric harmonicity of a whole tone (the interval between C and D) by $h_2$, and so on up to the symmetric harmonicity $h_{11}$ of a major seventh (the interval between C and B).

The basic philosophy of how we define $h_n$ follows Brefeld (though we omit the square roots since, as mentioned above, they make no difference for comparisons): we wish to define $h_n$ as the product of the numerator and denominator of the frequency ratio between the bottom note of the scale and the $n$th note.  However, there is a twist in our definition, in order to ensure that the inversion equivalence (b) will be satisfied when we pass to triads: inversions replace an interval with its additive inverse in $\ZZ_{12}$ (for example, the first F-major chord in Figure~\ref{fig:inversions} contains a perfect fourth between C and F, whereas the next inversion has a perfect fifth between F and the C above), so we will insist that $h_n = h_{12-n}$.  (This mirrors the fact that in $\ZZ_{12}$ there are only six meaningful nonzero distances.)  We will achieve this symmetry with a bit of brute force: by setting $h_n$ equal to whichever of the two possible choices for its value is smaller.

As an example, working again with Kepler's just scale from Table~\ref{tab:classicjustscale}, we see that the product of the numerator and denominator of the frequency ratio from the bottom of the scale to the fifth scale step (representing an interval of a perfect fourth) is $4 \cdot 3$, whereas from the bottom to the seventh scale step (representing a perfect fifth) it is $3 \cdot 2$.  We prefer the smaller of these options, so we set $h_5 = h_7 = 6$.

Putting it mathematically, we have arrived at the following definition.\footnote{This is far from the only reasonable way to assign a numerical value to the consonance of an interval.  For another example, see \cite{huron-intervalclass}, which assigns a numerical consonance to intervals by aggregating the results of several empirical studies \cite{malmberg, kameokakuriyagawa, hutchinsonknopoff}.}

\begin{definition}
Fix a just scale, and for each $n \in \ZZ_{12}$, denote by $\frac{a_n}{b_n} \in \QQ$ the frequency ratio between the $n$th note of the scale and the bottom ($0$th) note, in lowest terms.  For $1 \leq n \leq 11$, define
\[h_n = \min\left( a_n \cdot b_n, \; a_{12-n} \cdot b_{12-n}\right).\]
\end{definition}

Applying this to all of the intervals in Kepler's just scale, its symmetric harmonicities are shown in Table~\ref{tab:classicjustscaleintfractions}.
\setlength{\extrarowheight}{3pt}
\begin{table}[ht]
\begin{center}
\begin{tabular}{c|c|c|c|c|c|c|c|c|c|c}
$h_1$ & $h_2$ & $h_3$ & $h_4$ & $h_5$ & $h_6$ & $h_7$ & $h_8$ & $h_9$ & $h_{10}$ & $h_{11}$ \\
\hline 
$15 \cdot 8$ & $9 \cdot 8 $ & $5 \cdot 3 $ & $5 \cdot 4 $ & $3 \cdot 2 $ &$45 \cdot 32 $ &  $3 \cdot 2 $ & $5 \cdot 4 $ & $5 \cdot 3 $ & $9 \cdot 8 $ & $15 \cdot 8 $ 
\end{tabular}
\end{center}
\vspace{10pt}
\caption{The symmetric harmonicities stemming from Kepler's just scale.}\label{tab:classicjustscaleintfractions}
\end{table}

Next, we define the symmetric harmonicity of a triad to be the product of
the symmetric harmonicities of all the pairwise intervals between its notes. 
For example, the major triad C-E-G features a major third (four scale steps) from C to E, a perfect fifth (seven scale steps) from C to G, and a minor third (three scale steps) from E to G, so its symmetric harmonicity in Kepler's just scale equals $h_4 \, h_7 \, h_3 = 1800$.  Explicitly, the definition is the following.

\begin{definition}
For each $(a,b) \in T_0$ (which we view as a triad consisting of the $0$th, $a$th, and $b$th notes of the scale), the corresponding {\it symmetric harmonicity} is
\[h_{(a,b)} = h_a h_b h_{b-a}.\]
\end{definition}

By construction, this definition is invariant under all three types of equivalence described above.
Table~\ref{table:classicjustscaleprimecomplexity} shows the symmetric harmonicities (still based on Kepler's scale) for all twelve equivalence
classes: the first row shows a representative $(a,b) \in T_0$ of each equivalence class (so that
$a$ and $b$ are the distances from the base note to the other two notes in the triad), the second
row shows $h_{(a,b)}$ (in thousands, for simplicity), and the third row shows the resulting
ranking (so a smaller symmetric harmonicity indicates a better rank).
\setlength{\extrarowheight}{3pt}
\begin{table}[ht]
\begin{center}
\begin{tabular}{c|c|c|c|c|c|c|c|c|c|c|c}
$(1,2)$ & 
$(1,3)$ & 
$(1,4)$ & 
$(1,5)$ & 
$(1,6)$ & 
$(2,4)$ & 
$(2,5)$ & 
$(2,6)$ & 
$(2,7)$ & 
$(3,6)$ & 
$(3,7)$ & 
$(4,8)$ \\
\hline 
1036.8 & 129.6 & 36 & 14.4 & 1036.8 & 103.68 & 6.48 & 2073.6 & 2.592 & 324 & 1.8 & 8 \\
\hline
10 & 8 & 6 & 5 & 10 & 7 & 3 & 12 & 2 & 9 & 1 & 4
\end{tabular}
\end{center}
\vspace{10pt}
\caption{The symmetric harmonicities (in thousands) stemming from Kepler's just scale; the last row
gives the resulting ranking.}\label{table:classicjustscaleprimecomplexity}
\end{table}

\begin{figure}[ht] 
  \includegraphics[height=.52in]{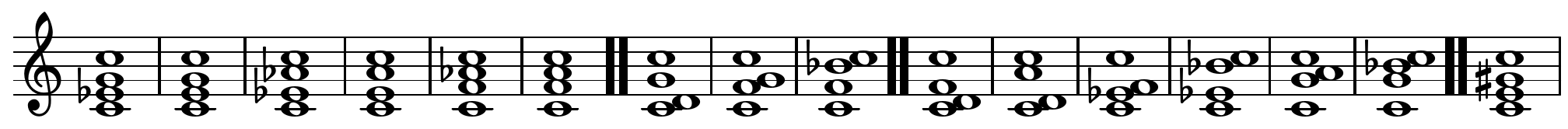}
\caption{The chords in the equivalence classes $(3,7)$, $(2,7)$, $(2,5)$, and $(4,8)$.}\label{fig:classicwinners}
\end{figure}

The symmetric harmonicity gives a (very rough) measure of the consonance/dissonance of a chord.
Indeed, the ``lightest'' chord family $(3,7)$ contains major and minor triads, which one might hope are the most consonant. The next three chord
families are $(2,7)$ (containing the suspended chords Csus2, Fsus2, and Fsus4), $(2,5)$ (another lovely chord family), and $(4,8)$ (corresponding to the augmented triad).
In Figure~\ref{fig:classicwinners}, we show all representatives based at C of each of these ``top'' equivalence classes.


\section{Variations on a Theme}

By design, the above analysis depends on the scale being used. To illustrate this, and to compare
the resulting ratings, we now collect data from various well-known just scales; for sake of comparison, we include also Kepler's just scale considered above.  The scales we consider are the following:
\vspace{12pt}
\begin{enumerate}[A:]
  \item Kepler's just scale (Table \ref{tab:classicjustscale}); \\
  \item Wendy Carlos's \emph{super just scale}~\cite[Section~6.1]{bensonmusic}, using primes $\le 17$; \\ 
  \item A continued fraction scale, approximating each power of $\sqrt[12] 2$ by the first convergent of its continued fraction that is within 1\% of the actual value;\\
  \item Stolzenburg's \emph{rational tuning}~\cite{stolzenburgjmm}, which is almost identical to
scale C,
with the only difference being in the semitone.\footnote{This difference does not affect symmetric harmonicities, which is why rows C and D of Table~\ref{tab:modharmvariousscales} are identical.} \\
  \item the Pythagorean scale with tritone frequency $729 \rat 512$, an ancient scale constructed via the circle of fifths.
\end{enumerate}
\vspace{12pt}
Their frequency ratios are listed in Table~\ref{tab:variousjustscales} and the resulting
symmetric harmonicities in Table~\ref{tab:modharmvariousscales}.
\begin{table}[ht]
\begin{center}
\begin{tabular}{c|c|c|c|c|c|c|c|c|c|c|c|c}
 & C & C\sh & D & E\fl & E & F & F\sh & G & G\sh & A & B\fl & B \\
\hline 
A & $1 \rat 1$ & $16 \rat 15$ & $9 \rat 8$ & $6 \rat 5$ & $5 \rat 4$ & $4 \rat 3$ & $45 \rat 32$ &
$3 \rat 2$ & $8 \rat 5$ & $5 \rat 3$ & $16 \rat 9$ & $15 \rat 8$ \\
\hline
B & $1 \rat 1$ & $17 \rat 16$ & $9 \rat 8$ & $6 \rat 5$ & $5 \rat 4$ & $4 \rat 3$ & $11 \rat
8$ & $3 \rat 2$ & $13 \rat 8$ & $5 \rat 3$ & $7 \rat 4$ & $15 \rat 8$ \\
\hline 
C & $1 \rat 1$ & $17 \rat 16$ & $9 \rat 8$ & $6 \rat 5$ & $5 \rat 4$ & $4 \rat 3$ & $17 \rat 12$ &
$3 \rat 2$ & $8 \rat 5$ & $5 \rat 3$ & $16 \rat 9$ & $15 \rat 8$ \\
\hline
D & $1 \rat 1$ & $16 \rat 15$ & $9 \rat 8$ & $6 \rat 5$ & $5 \rat 4$ & $4 \rat 3$ & $17 \rat 12$ &
$3 \rat 2$ & $8 \rat 5$ & $5 \rat 3$ & $16 \rat 9$ & $15 \rat 8$ \\
\hline
E & $1 \rat 1$ & $256 \rat 243$ & $9 \rat 8$ & $32 \rat 27$ & $81 \rat 64$ & $4 \rat 3$ & $729
\rat 512$ & $3 \rat 2$ & $128 \rat 81$ & $27 \rat 16$ & $16 \rat 9$ & $243 \rat 128$
\end{tabular}
\end{center}
\vspace{10pt}
\caption{The frequency ratios of the five scales.}\label{tab:variousjustscales}
\end{table}

\setlength{\extrarowheight}{3pt}
\begin{table}[ht]
\begin{center}
\begin{tabular}{c|c|c|c|c|c|c|c|c|c|c|c}
 & $h_1$ & $h_2$ & $h_3$ & $h_4$ & $h_5$ & $h_6$ & $h_7$ & $h_8$ & $h_9$ & $h_{10}$ & $h_{11}$ \\
\hline 
A & $15 \cdot 8$ & $9 \cdot 8 $ & $5 \cdot 3 $ & $5 \cdot 4 $ & $3 \cdot 2 $ &$45 \cdot 32 $ &  $3
\cdot 2 $ & $5 \cdot 4 $ & $5 \cdot 3 $ & $9 \cdot 8 $ & $15 \cdot 8 $ \\
\hline 
B & $15 \cdot 8$ & $7 \cdot 4 $ & $5 \cdot 3 $ & $5 \cdot 4 $ & $3 \cdot 2 $ & $11 \cdot 8 $ &  $3
\cdot 2 $ & $5 \cdot 4 $ & $5 \cdot 3 $ & $7 \cdot 4 $ & $15 \cdot 8 $ \\
\hline 
C & $15 \cdot 8$ & $9 \cdot 8 $ & $5 \cdot 3 $ & $5 \cdot 4 $ & $3 \cdot 2 $ &$17 \cdot 12 $ &  $3
\cdot 2 $ & $5 \cdot 4 $ & $5 \cdot 3 $ & $9 \cdot 8 $ & $15 \cdot 8 $ \\
\hline 
D & $15 \cdot 8$ & $9 \cdot 8 $ & $5 \cdot 3 $ & $5 \cdot 4 $ & $3 \cdot 2 $ &$17 \cdot 12 $ &  $3
\cdot 2 $ & $5 \cdot 4 $ & $5 \cdot 3 $ & $9 \cdot 8 $ & $15 \cdot 8 $ \\
\hline 
E & $243 \cdot 128$ & $9 \cdot 8 $ & $27 \cdot 16$ & $81 \cdot 64$ & $3 \cdot 2 $ & $729 \cdot
512$ &  $3 \cdot 2 $ & $81 \cdot 64 $ & $27 \cdot 16$ & $9 \cdot 8 $ & $243 \cdot 128$
\end{tabular}
\end{center}
\vspace{10pt}
\caption{The symmetric harmonicities stemming from the five scales.}\label{tab:modharmvariousscales}
\end{table}

We collect the triad ranking stemming from the symmetric harmonicity data for
the five scales in Table~\ref{table:variousscaleprimecomplexity}. Notice that there is a close agreement among the top three equivalence classes; beyond that, scales A--D continue to stay together, while the Pythagorean scale E features a noticeable different ranking,
stemming from the larger numbers present in its frequency ratios.
\setlength{\extrarowheight}{3pt}
\begin{table}[ht]
\begin{center}
\begin{tabular}{c|c|c|c|c|c|c|c|c|c|c|c|c}
&
$(1,2)$ & 
$(1,3)$ & 
$(1,4)$ & 
$(1,5)$ & 
$(1,6)$ & 
$(2,4)$ & 
$(2,5)$ & 
$(2,6)$ & 
$(2,7)$ & 
$(3,6)$ & 
$(3,7)$ & 
$(4,8)$ \\
\hline 
A & 10 & 8 & 6 & 5 & 10 & 7 & 3 & 12 & 2 & 9 & 1 & 4 \\
\hline 
B & 12 & 10 & 8 & 5 & 11 & 6 & 3 & 9 & 1 & 7 & 2 & 4 \\
\hline 
C & 12 & 9 & 6 & 5 & 10 & 8 & 3 & 11 & 2 & 7 & 1 & 4 \\
\hline 
D & 12 & 9 & 6 & 5 & 10 & 8 & 3 & 11 & 2 & 7 & 1 & 4 \\
\hline 
E & 7 & 5 & 7 & 5 & 7 & 4 & 2 & 11 & 1 & 7 & 3 & 11 
\end{tabular}
\end{center}
\vspace{10pt}
\caption{The symmetric harmonicity rankings stemming from the five scales.}\label{table:variousscaleprimecomplexity}
\end{table}


\section{Interlude: Comparison With Other Studies}

Our notion of symmetric harmonicity is clearly a cousin of Stolzenburg's harmonicity \cite{stolzenburgjmm}.  To compare the two directly,\footnote{We do not undertake a direct comparison between our theory and Brefeld's, as it is less well-known in the literature.} it is illuminating first to notice just how translation-dependent Stolzenburg's harmonicity values are.  As a concrete example, when one works with Kepler's just scale, we have computed previously that the major triad C-E-G has relative frequencies $\frac{1}{1}, \; \frac{5}{4}, \; \frac{3}{2}$, so its harmonicity equals
\[\text{lcm}(1,4,2)=4.\]
By contrast, the translated major triad E\fl-G\fl-B\fl\; has relative frequencies $\frac{6}{5}/\frac{6}{5}= \frac{1}{1}, \; \frac{45}{32}/\frac{6}{5} = \frac{75}{64}, \; \frac{16}{9}/\frac{6}{5} = \frac{40}{27}$, so its harmonicity equals
\[\text{lcm}(1,64,27) = 1728.\]
One could, of course, achieve a translation-invariant notion of harmonicity by averaging the harmonicities across all translations of the chord.  However, this average is increased dramatically by translations with high harmonicity, such as the above, and the result is an average that depends in quite an unpredictable way on the choice of scale.  To illustrate this, we have computed the average of the harmonicity values across all representatives of each of our equivalence classes---that is, across all translations and inversions of each triad---for the same scales considered above.  The following table records the rankings of each equivalence class according to these averages, where a ``1'' again indicates the lowest average (and therefore, in principle, the most ``pleasing'' triad as measured by average harmonicity).
\setlength{\extrarowheight}{3pt}
\begin{table}[ht]
\begin{center}
\begin{tabular}{c|c|c|c|c|c|c|c|c|c|c|c|c}
&
$(1,2)$ & 
$(1,3)$ & 
$(1,4)$ & 
$(1,5)$ & 
$(1,6)$ & 
$(2,4)$ & 
$(2,5)$ & 
$(2,6)$ & 
$(2,7)$ & 
$(3,6)$ & 
$(3,7)$ & 
$(4,8)$ \\
\hline
A &  5& 10  &4  & 8 & 7 & 12 & 1 & 6 &3  & 9 & 11 & 2\\
\hline
B &  9 & 10  & 3 & 5 & 8 & 12 & 6 & 11 & 4 & 7 & 2 &  1\\
\hline
C &  3 & 11  & 4 & 8 & 5 & 9 & 12 & 10 & 2 & 7 & 6 &  1\\
\hline
D &  9 & 12  & 6 & 8 & 5 & 7 & 3 & 11 & 1 & 10 & 4 &  2\\
\hline
E &  5 & 6  & 2 & 8 & 7 & 10 & 11 & 4 & 12 & 3 & 9 &  1\\
\end{tabular}
\end{center}
\vspace{10pt}
\caption{The rankings of Stolzenburg harmonicities (averaged across all representatives of each equivalence class) stemming from various scales.}
\label{tab:averageStolzenburg}
\end{table}

These rankings seem significantly more sensitive to the choice of scale than do those of Table~\ref{table:variousscaleprimecomplexity}.  Furthermore, it is intriguing to notice that the augmented triad (4,8) outperforms the major/minor triad (3,7) in every scale by this metric---often quite dramatically, such as in the Kepler just scale, where the major/minor triad ranks near the bottom.  Heuristically, this is due to the fact that the (4,8) equivalence class has the fewest representatives, so it is the least susceptible to having its average harmonicity affected by particularly ``bad'' representatives of the equivalence class.

The surprising behavior in Table~\ref{tab:averageStolzenburg} should not be taken as a sign of deficiency in Stolzenburg's notion of harmonicity, but rather, as an indication that averaging across translations is perhaps not the right way to achieve translation-invariance after all; instead (as we have done in our notion of symmetric harmonicity), we advocate for imposing translation-invariance by fiat, simply by treating every instance of, for example, a major third as the same interval.  In fact, a similar idea seems to be implicit in Stolzenburg's work as well; for example, in \cite[Table 4]{stolzenburgjmm}, a major triad is indicated with scale steps $\{0,4,7\}$ regardless of its base note.

In addition to Stolzenburg's work on harmonicity, another study to which we should compare our work is the empirical data on how pleasant different chords are actually perceived to be.   The most widely-referenced empirical study of harmony
perception of musical chords seems to be Johnson-Laird, Kang, and Leong's~\cite{johnsonlairdkangleong}. We cannot compare our
rankings directly to those in that study, since the ratings
of~\cite[Experiment~1]{johnsonlairdkangleong} differ among inversions of the same chord. (To make
things even more subtle, each inversion features a different base note, and the rankings differ
quite prominently among inversions of the same chord. This makes for a fascinating study involving human subjects, but it is somewhat orthogonal to our approach/assumptions.)  To account for this difference, similarly to the above, we have computed the average rating according to \cite{johnsonlairdkangleong} among all representatives within a given equivalence class that were included in their experiment.  The results are recorded in Table~\ref{tab:empirical}.  This study used an (approximately) equal-tempered scale, so for comparison, we also reproduce the symmetric harmonicity rankings in the continued fraction scale (the rational scale that is ``closest'' to equal-tempered) in the second row of the table.

In the third row of Table~\ref{tab:empirical}, we have also computed the rankings of the triads in \cite[Experiment~1]{johnsonlairdkangleong} according to Stolzenburg's harmonicity.  In light of the above discussion, however, we have computed these in a translation-invariant way: for each chord appearing in \cite[Experiment~1]{johnsonlairdkangleong}, we have computed its harmonicity in a scale whose base note is taken to be the lowest note of the chord.  Because the empirical study of Johnson-Laird et al includes multiple representatives of each of our equivalence classes, we have then averaged the harmonicity values across all representatives appearing in the study to obtain the rankings in Table~\ref{tab:empirical}.

\setlength{\extrarowheight}{3pt}
\begin{table}[ht]
\begin{center}
\begin{tabular}{c|c|c|c|c|c|c|c|c|c|c|c|c}
&
$(1,2)$ & 
$(1,3)$ & 
$(1,4)$ & 
$(1,5)$ & 
$(1,6)$ & 
$(2,4)$ & 
$(2,5)$ & 
$(2,6)$ & 
$(2,7)$ & 
$(3,6)$ & 
$(3,7)$ & 
$(4,8)$ \\
\hline 
empirical & 10 & 6 & n/a & 7 & 8 & n/a & 3 & 5 & 2 & 4 & 1 & 9\\
\hline
symm harm & 12 & 9 & 6 & 5 & 10 & 8 & 3 & 11 & 2 & 7 & 1 & 4 \\
\hline
avg harm &  9 & 10  & n/a & 8 & 6 & n/a & 4 & 7 & 2 & 4 & 1 &  3\\
\end{tabular}
\end{center}
\vspace{10pt}
\caption{The empirical data, as compared to the symmetric harmonicity rankings of our work and the Stolzenburg harmonicity ratings (both in the continued fraction scale).  Harmonicity ratings are computed in a scale with base note equal to the lowest note in the chord, and are averaged when the empirical study includes multiple representatives of the same equivalence class.}
\label{tab:empirical}
\end{table}

Several features of this table are worth noting.  First, both symmetric harmonicity and harmonicity do a fairly good job of reproducing the empirical rankings, especially for the highest-ranked chords; this is certainly a valuable sanity check.  Furthermore, the surprising behavior of Stolzenburg's harmonicity noted in Table~\ref{tab:averageStolzenburg} disappears when we avoid averaging over different translations.  This, to us, further emphasizes the advantage of building translation-invariance (as well as the other types of invariance discussed in this work) directly into the theory, rather than relying on averaging to achieve it.


\section{Coda: 4-Chords}

Chords with four distinct notes are prominent in many musical genres, and as we mentioned above,
our methodology applies to them just as well as to triads. There are now many more equivalence
classes, and so this section gives only a glimpse.

Table~\ref{tab:symmharm4winners} shows the symmetric harmonicities of all equivalence classes of
4-chords for which this quantity is below $10^8$. They are based on either scale A, C, or D
(their symmetric harmonicities differ only for the tritone 6, which does not appear as a
difference in the top-7 list).
\setlength{\extrarowheight}{3pt}
\begin{table}[ht]
\begin{center}
\begin{tabular}{c|c|c|c|c|c|c}
(2,5,9) & (2,5,7) & (1,5,8) & (2,4,7) & (1,4,9) & (1,4,8) & (1,3,8) \\
\hline 
11.7 & 16.8 & 25.9 & 56.0 & 64.8 & 86.4 & 93.3
\end{tabular}
\end{center}
\vspace{10pt}
\caption{The symmetric harmonicities of some 4-chords, in millions.}\label{tab:symmharm4winners}
\end{table}

\begin{figure}[ht] 
  \includegraphics[height=.53in]{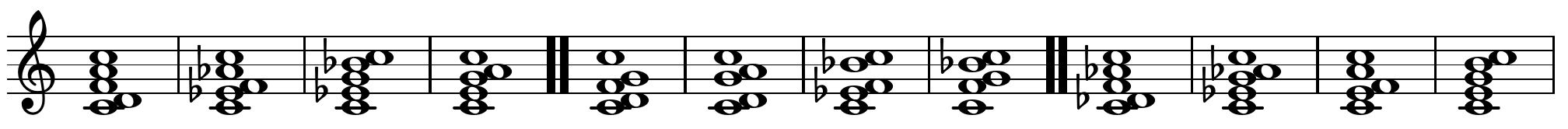}
\caption{The chords in the equivalence classes $(2,5,9)$, $(2,5,7)$, and
$(1,5,8)$.}\label{fig:4chordwinners}
\end{figure}

The ``winning'' equivalence class contains both the minor 7th and major 6th chords; the third
class contains the major 7th chord. The top three chords are depicted in Figure~\ref{fig:4chordwinners}.

The first four chord families occupy ranks 4, 6, 2, and 1, respectively,
in the empirical study~\cite[Experiment~2]{johnsonlairdkangleong}, which is somewhat suggestive that symmetric harmonicity is in some agreement with perceptions even for 4-chords.  A notable absence in Table~\ref{tab:symmharm4winners}, on the other hand, is the dominant 7th chord (3,5,9) with a symmetric harmonicity of 397 millions, which ranked fifth in the empirical study.  We end our exploration here, but we hope that this glimpse of uncharted territory serves as an invitation for further investigation.


\bibliographystyle{amsplain}
\bibliography{bib}

\setlength{\parskip}{0cm} 

\end{document}